\documentclass[10pt,a4paper]{amsart}
\usepackage[T1]{fontenc}
\usepackage{%
	lmodern,%
	mparhack,%
	relsize%
}
\usepackage[english]{babel}
\usepackage[autostyle=true]{csquotes}
\usepackage{%
	amsaddr,%
	amsfonts,%
	amsmath,%
	amsthm,%
	amscd,%
	amssymb,%
	mathtools%
}
\usepackage{amstext}
\usepackage{etoolbox}
\usepackage[ruled]{algorithm2e}
\usepackage{float}
\usepackage{tabularx}
\usepackage{xcolor}   			
\usepackage{graphicx}
\graphicspath{ {./images/} }
\usepackage[hyphens]{url}
\usepackage{hyperref}
\hypersetup{%
	naturalnames=false,%
	breaklinks=true,%
	colorlinks=true,
	frenchlinks=false,%
	linkcolor=red,%
	citecolor=red,%
	debug=true,%
	filecolor=magenta,%
	urlcolor=magenta,%
	linktoc=all,
}
\usepackage{orcidlink} 			
\usepackage[%
	backend=biber,%
	citestyle=numeric-comp,%
	style=numeric-comp,%
	sorting=nty,%
	hyperref=true,%
	block=space,%
	backref=false%
]{biblatex}
\usepackage{marvosym}

\newcolumntype{Z}{>{\setbox0=\hbox\bgroup}c<{\egroup}@{\hspace*{-\tabcolsep}}} 

\pretocmd\mvchr{\text}{}{\errmessage{Patching \noexpand\mvchr failed}}
\pretocmd\textmvs{\text}{}{\errmessage{Patching \noexpand\textmvs failed}}

\newcommand{\algref}[1]{\hyperref[alg:#1]{algorithm~\ref{alg:#1}}}
\newcommand{\chapref}[1]{\hyperref[chap:#1]{chapter~\ref{chap:#1}}}
\newcommand{\secref}[1]{\hyperref[sec:#1]{section~\ref{sec:#1}}}
\newcommand{\subsecref}[1]{\hyperref[subsec:#1]{subsec.~\ref{subsec:#1}}}
\newcommand{\subsubsecref}[1]{\hyperref[subsubsec:#1]{subsubsec.~\ref{subsubsec:#1}}}
\newcommand{\figref}[1]{\hyperref[fig:#1]{fig.~\ref{fig:#1}}}
\newcommand{\tabref}[1]{\hyperref[tab:#1]{table~\ref{tab:#1}}}
\newcommand{\defref}[1]{\hyperref[def:#1]{definition~\ref{def:#1}}}
\newcommand{\thmref}[1]{\hyperref[thm:#1]{theorem~\ref{thm:#1}}}
\newcommand{\lemmaref}[1]{\hyperref[lemma:#1]{lemma~\ref{lemma:#1}}}
\newcommand{\corollaryref}[1]{\hyperref[corollary:#1]{corollary~\ref{corollary:#1}}}
\newcommand{\remarkref}[1]{\hyperref[rmk:#1]{remark~\ref{rmk:#1}}}
\newcommand{\exampleref}[1]{\hyperref[example:#1]{example~\ref{example:#1}}}

\makeatletter
\newcommand{\homcoords}[1]{\left[#1\checknextarg}
\newcommand{\checknextarg}{\@ifnextchar\bgroup{\gobblenextarg}{\right]}}
\newcommand{\gobblenextarg}[1]{ : #1\@ifnextchar\bgroup{\gobblenextarg}{\right]}}
\makeatother

\newcommand{\commandof}[2]{\mathmbox{#1\roundbrack{#2}}}
\newcommand{\quotient}[2]{\mathmbox{#1 / #2}}
\newcommand{\intring}[1]{\mathmbox{\mathcal{O}_{#1}}}

\newcommand{\primeidlocal}[1]{\mathmbox{\mathfrak{m}_{#1}}}
\newcommand{\uniformizerlocal}[1]{\mathmbox{\varpi_{#1}}}
\newcommand{\residuefield}[1]{\quotient{\intring{#1}}{\primeidlocal{#1}}}

\newcommand{\brackets}[1]{\mathmbox{\left\lbrace#1\right\rbrace}}

\newcommand{\roundbrack}[1]{\mathmbox{\left(#1\right)}}

\newcommand{\modp}[1]{\ \mathmbox{\roundbrack{\mathrm{mod}\ #1}}}
\newcommand{\ordp}[1]{\mathrm{ord}\roundbrack{#1}}
\newcommand{\floor}[1]{\left\lfloor#1\right\rfloor}
\newcommand{\ceil}[1]{\left\lceil#1\right\rceil}

\newcommand{\Ker}{\mathrm{Ker}}
\newcommand{\Kerof}[1]{\commandof{\Ker}{#1}}
\newcommand{\Img}{\mathrm{Im}}
\newcommand{\Imgof}[1]{\commandof{\Img}{#1}}
\newcommand{\Mod}{\mathrm{Mod}}
\newcommand{\Modof}[1]{\commandof{\Mod}{#1}}
\newcommand{\QMod}{\mathrm{Mod}}
\newcommand{\QModof}[2]{\commandof{\QMod_{#1}}{#2}}
\newcommand{\QModx}[1]{\QMod_{p^{#1}}}
\newcommand{\QModxof}[2]{\commandof{\QModx{#1}}{#2}}
\newcommand{\Exp}{\mathrm{Exp}}
\newcommand{\Expof}[1]{\commandof{\Exp}{#1}}
\newcommand{\Expinv}{{\Exp}^{-1}}
\newcommand{\Expinvof}[1]{\commandof{\Expinv}{#1}}
\newcommand{\Log}{\mathrm{Log}}
\newcommand{\Logof}[1]{\commandof{\Log}{#1}}
\newcommand{\degof}[1]{\commandof{\mathrm{deg}}{#1}}
\newcommand{\charof}[1]{\commandof{\mathrm{char}}{#1}}
\newcommand{\cardof}[1]{\commandof{\mathrm{card}}{#1}}
\newcommand{\C}{\mathbb{C}}

\newcommand{\Q}{\mathbb{Q}}
\newcommand{\F}{\mathbb{F}}
\newcommand{\Z}{\mathbb{Z}}
\newcommand{\N}{\mathbb{N}}
\newcommand{\K}{\mathbb{K}}
\newcommand{\Zof}[1]{\quotient{\Z}{#1 \Z}}

\newcommand{\weier}{\mathcal{W}}
\newcommand{\hyperec}{\mathcal{H}}
\newcommand{\jacobian}{\mathcal{J}}
\newcommand{\jacobianof}[1]{\commandof{\jacobian}{#1}}
\newcommand{\qpjacobianof}[2]{\commandof{\jacobian_{#1}}{#2}}
\newcommand{\fullcurveof}[2]{\commandof{#1}{#2}}
\newcommand{\weierof}[1]{\fullcurveof{\weier}{#1}}
\newcommand{\hyperecof}[1]{\fullcurveof{\hyperec}{#1}}
\newcommand{\suppof}[1]{\commandof{\mathrm{supp}}{#1}}
\newcommand{\Qp}[1]{\Q_{#1}}
\newcommand{\GF}[1]{\commandof{\mathrm{GF}}{#1}}

\newtheorem{theorem}{Theorem}
\numberwithin{theorem}{section} 

\newtheorem{remark}{Remark}
\numberwithin{remark}{section} 

\newtheorem{definition}{Definition}
\numberwithin{definition}{section} 
\newtheorem{example}{Example}

\numberwithin{equation}{section}
\numberwithin{figure}{section}

\title[On the DLP for elliptic curves over local fields]{%
	On the Discrete Logarithm Problem for elliptic curves over local fields}
\author{Giuseppe Filippone \orcidlink{0000-0001-7315-1852} \href{mailto:giuseppe.filippone01@unipa.it}{$ \Letter $}}
\address{%
	\normalsize Dipartimento di Matematica e Informatica,\\%
	\normalsize Università degli Studi di Palermo\\%
	\normalsize Via Archirafi 34, 90123 Palermo, Italy}
\date{}

\begin{document}
	\begin{abstract}
		The Discrete Logarithm Problem (DLP) for elliptic curves has been extensively studied
		since, for instance, it is the core of the security of cryptosystems like Elliptic Curve
		Cryptography (ECC). In this paper, we present an attack to the DLP for elliptic curves
		based on its connection to the problem of lifting,
		by using the exponential map for elliptic curves and its inverse over $ \Zof{p^k} $.
		Additionally, we show that hyperelliptic curves are resistant to this attack, meaning that these
		latter curves offer a higher level of security compared to the classic elliptic curves
		used in cryptography.
	\end{abstract}

	\maketitle


	\section{Introduction}\label{sec:intro}

	One of the key areas of cryptography is the Elliptic Curve Cryptography (ECC),
	which has been widely adopted in many real-world applications due to its security
	and efficiency properties, since it requires smaller key sizes than other
	public-key cryptosystems, such as RSA. This makes ECC particularly
	well-suited for resource-constrained devices, such as smart cards, mobile phones,
	and embedded systems.
	The security of elliptic curve cryptography is based on the intractability of the
	Elliptic Curve Discrete Logarithm Problem (ECDLP) which,
	together with the integer factorization problem, is one of the fundamental problems
	underpinning the security of modern cryptography.
	Specifically, the ECDLP involves computing the scalar $ a \in \K $ such that
	$ Q = a \cdot P $, where $ P $ and $ Q $ are (publicly shared) points
	belonging to an elliptic curve defined over a field $ \K $.
	This problem is believed to be intractable, meaning that there is no known
	algorithm that can solve it efficiently, which makes Elliptic Curve Cryptography (ECC)
	a secure and reliable approach for secure communications.
	Specifically, the basic idea behind ECC is to use the difficulty of solving the DLP
	on elliptic curves to create a secure and efficient cryptographic system.

	However, despite the fact that the DLP on elliptic curves is widely believed
	to be intractable, there are certain methods that can be used to solve it
	efficiently for some classes of curves, i.e. curve with a large cofactor,
	a small cardinality, a small prime field characteristic, and so on.
	For instance, the ECDLP is easier in the case of anomalous elliptic curves,
	i.e. such that the trace of the Frobenius endomorphism of elliptic curve
	is equal to one \cite{Smart1999, LEPREVOST2005225, SamaevAnomalous, SatohArakiAnomalous}.
	Additionally, the Index Calculus algorithm and the Pollard's rho algorithm
	can be used to find a solution to the ECDLP on elliptic curves with small characteristic.
	Nevertheless, beside the case of anomalous elliptic curves, these methods are
	not practical for most cases of practical interest, since they require a large amount
	of memory and computation, and their complexity grows exponentially with the
	number of points belonging to the elliptic curve.

	It is known that the Weierstrass $ \wp $-function and its derivative $ \wp^\prime $,
	over the field of complex numbers, define an exponential map
	$ \Exp \colon z \mapsto \roundbrack{\wp(z), \frac{1}{2} \wp^\prime(z)} $
	(also known as Weierstrass uniformizing map)
	between the elements $ z $ belonging to the complex torus
	$ \C / \Lambda $, where $ \Lambda $ is the period lattice of $ \wp $,
	and the points belonging to an elliptic curve in short Weierstrass form, given
	by the equation $ {\roundbrack{\frac{1}{2} \wp^\prime(z)}}^2 = \wp^3(z) - \frac{g_2}{4} \wp(z) - \frac{g_3}{4} $,
	where $ g_2, g_3 \in \C $ are known as elliptic invariants,
	of the corresponding complex projective plane.
	This map has the following property (see e.g. $ \S $VI and $ \S $IX in \cite{Silverman2009}),
	i.e. $ \Expof{z_1 + z_2} = \Expof{z_1} \mathbin{\star} \Expof{z_2} $,
	where the operation $ (\star) $	is the addition on the points of
	the elliptic curve by using the well-known chord-and-tangent law.

	An authoritative literature on the connection between the problem of lifting
	and the ECDLP, summarized in \cite{SilvermanLiftingECDLP}, motivated us
	to provide, in this paper, a cryptanalytic framework to attack
	the DLP for elliptic curves in Weierstrass form
	over local fields by using the map $ \Exp $,
	and we confine ourselves to the local field $ \Qp{p} $ of $ p $-adic numbers.
	Additionally, we show that, as a consequence of the fact the greatest common divisor (GCD) and
	the reduction modulo $ p^k $ do not commute over $ \K[x] $, where $ \K $ is a field,
	our framework is not applicable to hyperelliptic curves, meaning that these latter curves
	are resistant to these kinds of attacks.

	Furthermore, cryptosystems over local fields have recently received attention.
	For instance, in \cite{Xu2008} the authors gave a cryptosystem based on quotient
	groups of an elliptic curve over $ \Qp{p} $, encrypting messages with
	variable lengths, which led to public-key cryptosystems with hierarchy management
	\cite{XuYue}.

	\section{Preliminaries}\label{sec:notations}

	In this section, we provide to the reader the notation used throughout this
	paper, and the methodologies used deal with the ECDLP. We address the reader
	to classic book, e.g. \cite{Husemoller, Mumford, Silverman2009},
	for further details about elliptic curves.

	\begin{definition}[{\normalfont cf. \cite{Silverman2009}}]
		Let $ \K $ be a field. An elliptic curve $ \weier $ in (generalized) Weierstrass form
		is a (smooth) algebraic curve over $ \K $ of genus $ g = 1 $ given by the equation
		\begin{equation}\label{eq:general-ec}
			y^2 + a_1 x y + a_3 y = x^3 + a_2 x^2 + a_4 x + a_6.
		\end{equation}
		Denote by $ \Omega = \homcoords{Z}{X}{Y} = \homcoords{0}{0}{1} $ its
		point at infinity.
	\end{definition}

	If $ \charof{\K} \ne 2, 3 $, then the map
	\begin{equation}\label{eq:tranform-map}
		 (x, y) \mapsto \roundbrack{f(x) - {\roundbrack{\frac{h(x)}{2}}}^2, y - \frac{h(x)}{2}},
	\end{equation}
	where $ f(x) = x^3 + a_2 x^2 + a_4 x + a_6 $, and $ h(x) = a_1 x + a_3 $,
	transforms the equation \eqref{eq:general-ec} to the form
	\begin{equation}\label{eq:short-ec}
		y^2 = x^3 + a x + b,
	\end{equation}
	known as short Weierstrass form of an elliptic curve.

	Let $ \K $ be a local field, $ \intring{\K} $ its ring of integers,
	$ \primeidlocal{\K} $ its prime ideal, and $ \mathsf{k} = \residuefield{\K} $
	its residue field (see, e.g. \cite{SerreLocalField} for a general reference).
	We recall that, given an elliptic curve $ \weier $ in short Weierstrass form,
	given by the equation \eqref{eq:short-ec}, whose reduction modulo
	$ \primeidlocal{\K} $ is non-singular, the following short sequence:
	\begin{equation}\label{eq:exact-sequence-w}
		0 \longrightarrow \primeidlocal{\K} \xrightarrow{\Exp} \qpjacobianof{\K}{\weier} \ \xrightarrow{\Mod} \qpjacobianof{\mathsf{k}}{\weier} \longrightarrow 0
	\end{equation}
	is exact \cite{kosterspannekoek, Hiranouchi2017LocalTP, rosawinter} (see also \cite[ch. $ \S $VII]{Silverman2009}),
	where $ \qpjacobianof{\K}{\weier} $ is the Jacobian of $ \weier $ over $ \K $ and
	$ \qpjacobianof{\mathsf{k}}{\weier} $ is the Jacobian of $ \weier $ over $ \mathsf{k} $.
	Hence, $ \Imgof{\Exp} = \Kerof{\Mod} $,
	$ \Exp $ is a monomorphism, $ \Mod $ is an epimorphism, and one has that
	\begin{equation}
		\qpjacobianof{\mathsf{k}}{\weier} \cong \quotient{\qpjacobianof{\K}{\weier}}{\Kerof{\Mod}} = 	\quotient{\qpjacobianof{\K}{\weier}}{\Imgof{\Exp}}.
	\end{equation}

	The map $ \Mod $, which is surjective by Hensel's lemma
	(see the proof in section $ \S $VII.2.1 of \cite{Silverman2009}),
	is the reduction modulo $ \primeidlocal{\K} $ of the coordinates
	of the points $ P = \homcoords{Z}{X}{Y} $ in $ \weierof{\K} $ which, up to a multiplication times a suitable
	element belonging to $ \intring{\K} $, have integral entries $ Z, X, Y \in \intring{\K} $:
	\begin{alignat}{2}\label{eq:ec-mod-map}
		\Mod \colon \ &\qpjacobianof{\K}{\weier}&\ \rightarrow\ &\qpjacobianof{\mathsf{k}}{\weier}\nonumber\\
		&\homcoords{Z}{X}{Y} - \Omega &\ \mapsto\ &\homcoords{Z \modp{\primeidlocal{\K}}}{X \modp{\primeidlocal{\K}}}{Y \modp{\primeidlocal{\K}}} - \Omega.
	\end{alignat}

	As we mentioned in the introduction, over the field $ \C $ of complex numbers,
	the function $ \Exp $ is defined as
	$ \Exp \colon z \mapsto \homcoords{1}{\wp(z)}{\frac{1}{2} {\wp^\prime(z)}} $.
	Moreover, there are the Laurent series expansions for $ \wp $ and
	$ \wp^\prime $ \cite{AbramowitzStegun1972, ApostolTomMike1990MfaD} given by
	\begin{align}\label{eq:wp-laurent}
		\wp(z) &= \frac{1}{z^2} + \sum_{l=2}^{\infty} c_l z^{2l - 2},\\\nonumber
		\wp^\prime(z) &= -\frac{2}{z^3} + \sum_{l=2}^{\infty} (2l - 2) c_l z^{2l - 3},
	\end{align}
	where $ c_2 = - \frac{1}{5} a $, $ c_3 = - \frac{1}{7} b $, and
	$ c_l = \frac{3}{(2l + 1)(l - 3)} \sum_{s=2}^{l - 2} c_s c_{l - s} $, which
	converge locally uniformly absolutely in $ \C \setminus \Lambda $.
	Similarly, the same series converge over $ \K $ for $ z $ belonging the maximal ideal
	$ \primeidlocal{\K} $, which plays here the role of a neighborhood of zero.
	Hence, one may define, over $ \K $, the homomorphism $ \Exp $ as
	\begin{alignat}{2}\label{eq:ec-exp-map}
		\Exp \colon \ &\primeidlocal{\K} &\ \longrightarrow\ &\qpjacobianof{\K}{\weier}\nonumber\\
		&z &\ \longmapsto\ &\homcoords{1}{\wp(z)}{\frac{1}{2} {\wp^\prime(z)}} - \Omega, \quad \text{if } z \ne 0,\\
		&0 &\ \longmapsto\ &\Omega - \Omega \nonumber.
	\end{alignat}

	Since $ z = 0 $ is the only element of $ \primeidlocal{\K} $
	mapped to $ \Omega - \Omega $, the homomorphism $ \Exp $
	is into and, for any $ z $ in the neighborhood of zero, one can define
	(see $ \S $IV and $ \S $VII \cite{Silverman2009} for further details)
	the function $ \Expinv $ such that, for $ z \in \primeidlocal{\K} $, gives in turn
	$ \Expinvof{\Expof{z}} = z $.
	In particular, for $ 1 < i \le 5 $, one has that the function
	\begin{equation}\label{eq:inv-exp}
		\Expinv \colon \homcoords{T}{X}{Y} \longmapsto -2\; \frac{X}{Y}
	\end{equation}
	is such that $ \Expinvof{\Expof{z}} = z \modp{{\primeidlocal{\K}}^i} $.
	As long as the domain of $ \Expinv $ is $ \Imgof{\Exp} $,
	one has that $ -2 \frac{X}{Y} $	is equivalent to $ -2 \frac{\wp(z)}{\wp^\prime(z)} $,
	whose series expansion is
	\begin{equation}\label{eq:series-wp-wp-der}
		z - \frac{2}{5} a z^5 - \frac{3}{7} b z^7 + O(z^9).
	\end{equation}
	Hence, for any $ z $ such that $ \Expof{z} = P = \homcoords{T}{X}{Y} $
	one has that $ z - \roundbrack{- 2 \frac{X}{Y}} \in {\primeidlocal{\K}}^i $,
	for $ 1 < i \le 5 $.


	It is worth to note that the exact sequence in equation \eqref{eq:exact-sequence-w}
	does not split over $ \mathsf{k} $ if one supposes that the elliptic curve
	in Weierstrass form is an anomalous curve
	(see example 1 in \cite{filippone2023expEA}). In particular, it holds the
	following theorem.

	\begin{theorem}[\normalfont cf. \cite{filippone2023expEA}]\label{thm:jw-splitting}
		If $ \mathsf{k} = \residuefield{\K} $ is finite, and $ \weier $ is not an anomalous curve,
		then $ \qpjacobianof{\K}{\weier} $ is isomorphic to the direct sum of
		$ \qpjacobianof{\mathsf{k}}{\weier} $ and $ \primeidlocal{\K} $.
	\end{theorem}

	\section{The exponential map over \texorpdfstring{$ \Qp{p}$}{Qp}}\label{sec:ell-qp}

	In this section, we explicit the maps $ \Mod $ and $ \Exp $
	given in equations \eqref{eq:ec-mod-map} and \eqref{eq:ec-exp-map}, respectively,
	in the case in which $ \K $ is the local field $ \Qp{p} $ of $ p $-adic
	numbers. In particular, from now on, we confine ourselves to this latter case.

	In the case in which $ \K = \Qp{p} $, one has that its ring of integers $ \intring{\K} $
	is the ring $ \Z_p $ of $ p $-adic integers, the prime ideal $ \primeidlocal{\K} $
	is $ p \Z_p $ (which uniformizer $ \uniformizerlocal{\K} $ is equal to $ p $), and the residue field
	$ \residuefield{\K} $ is $ \Zof{p} = \GF{p} $.

	Specifically, in this section, we study the field $ \Qp{p} $ and the sequence
	in equation \eqref{eq:exact-sequence-w} by using an inverse limit process,
	i.e. we approach to $ \Z_p $ starting from $ \Zof{p} $, and ending with
	$ \Zof{p^k} $, with $ k \rightarrow \infty $.

	We remark that going from $ \Zof{p} $ to $ \Zof{p^k} $ involves moving from a
	field to a ring in which there are elements that may not have an inverse.
	Therefore, rather than working with points in non-homogeneous coordinates $ (x, y) $,
	we will consider both the elliptic curve and the points in homogeneous coordinates, i.e.
	$ P = \homcoords{Z}{X}{Y} \equiv \roundbrack{\frac{X}{Z}, \frac{Y}{Z}} = (x, y) $.

	Since, in this case, the local field is $ \Qp{p} $, a neighborhood of zero
	is given by multiples of $ \uniformizerlocal{\Qp{p}} = p $, that is, any
	$ z $ such that $ p \mid z $. Hence, the series in equation \eqref{eq:wp-laurent}
	always converge modulo $ p^k $ since, in this neighborhood,
	$ c_l z^{2l - 2} \equiv 0 \modp{p^k} $, and $ (2l - 2) c_l z^{2l - 3} \equiv 0 \modp{p^k} $,
	for a suitable positive integer $ l $.

	Thus, the map $ \Exp $ is given by
	\begin{alignat}{2}\label{eq:exp-qp-w}
		\Exp \colon \ &\quotient{p \Z_p}{p^k \Z_p} &\ \longrightarrow\ &\weierof{\Zof{p^k}}\nonumber\\
		&z = ph &\ \longmapsto\ & \homcoords{1}{\wp(z)}{\frac{1}{2} {\wp^\prime(z)}} - \Omega,
	\end{alignat}
	where $ h = 1, 2, \ldots, p^{k - 1} $, since $ \Z_p = \lim\limits_{\longleftarrow} \Zof{p^k} $.
	Note that in this case we are considering the quotient $ \quotient{p \Z_p}{p^k \Z_p} $, since
	$ \Expof{ph} = \Expof{p h + p^k} $.
	Additionally, since $ \quotient{p \Z_p}{p^k \Z_p} $ is isomorphic to $ \Zof{p^{k - 1}} $
	through the mapping $ p h \mapsto h $, we will consider this latter as the domain
	of $ \Exp $.

	By multiplying the homogeneous coordinates of the points belonging to $ \Imgof{\Exp} $
	by a suitable integer such that each coordinate is integral, one notes that
	$ P \in \Imgof{\Exp} $ is such that $ P = \homcoords{p h_1}{p h_2}{h_3} $, where
	$ h_1, h_2, h_3 \in \brackets{1, \ldots, p^{k - 1}} $, and $ p \nmid h_3 $.
	Moreover, since $ \Expof{z_1} \not\equiv \Expof{z_2} \modp{p^k} $
	if $ z_1 \not\equiv z_2 \modp{p^k} $,
	the cardinality of $ \Imgof{\Exp} $ is $ p^{k - 1} $. Hence, the cardinality
	of $ \qpjacobianof{\Zof{p^k}}{\weier} $ is equal to
	$ \cardof{\qpjacobianof{\Zof{p}}{\weier}} \cdot p^{k - 1} $.

	Finally, we define the following function
	\begin{alignat}{2}
		\QModx{k} \colon \ &\qpjacobianof{\K}{\weier}&\ \longrightarrow\ &\qpjacobianof{\Zof{p^k}}{\weier}\nonumber\\
		&\homcoords{Z}{X}{Y} - \Omega &\ \longmapsto\ &\homcoords{Z \modp{p^k}}{X \modp{p^k}}{Y \modp{p^k}} - \Omega.
	\end{alignat}

	\section{The \texorpdfstring{$\Log$}{Log} function}\label{sec:ell-log}

	In this section, we provide an explicit map to compute the map $ \Expinv $ in
	equation \eqref{eq:inv-exp}, here denoted as $ \Log $, over $ \Zof{p^k} $
	in order to define our attack based on ECDLP in \secref{ell-ecdlp}. In particular,
	it holds the \thmref{log-map}.

	\begin{remark}
		We remark that, over $ \Zof{p^k} $, the map in equation \eqref{eq:inv-exp}
		do not give in turn the correct $ z $, but an element which is equivalent to it
		(at most) modulo $ p^5 $, since the series expansion of $ -2 \frac{\wp(z)}{\wp^\prime(z)} $
		is equal to $ z $ modulo $ p^i $, for $ 1 < i \le 5 $.
	\end{remark}

	\begin{remark}
		We remark that, although there is the inverse $ \wp^{-1} $ of the Weierstrass
		$ \wp $-function, the Puiseux series expansion of $ \wp^{-1}(z) $ involves square roots
		of powers of $ z $, which are not well defined over $ \Zof{p^k} $. In particular,
		for $ p \mid a $, it holds that the two operations $ {\roundbrack{a^{\frac{1}{c}}}}^b $ and
		$ {\roundbrack{a^b}}^{\frac{1}{c}} $ could give in turn two different results.
		For instance, consider $ p = 3 $, $ k = 4 $, and the operation $ 9^{\frac{5}{2}} $.
		On one hand, $ {\roundbrack{9^{\frac{1}{2}}}}^5 $ is equal to $ 0 $ modulo $ 3^4 $.
		On the other hand, $ {\roundbrack{9^5}}^{\frac{1}{2}} $ is equal to the square roots of
		$ 0 $ modulo $ 3^4 $, i.e. $ \brackets{0, 9, 18, 27, 36, 45, 54, 63, 72} $.
	\end{remark}

	\begin{theorem}[$ \Log $ map]\label{thm:log-map}
		The map $ \Exp $, given in equation \eqref{eq:exp-qp-w}, over $ \Zof{p^k} $
		is (locally) invertible for any neighborhood of zero given by any $ z \in \Z $
		such that $ p \mid z $. Denote this inverse map as $ \Log $.
	\end{theorem}
	\begin{proof}
		In this proof, we will show that it is possible to compute the map $ \Log $ by first determining
		the value $ z $ modulo $ p^5 $ from the series expansions of $ \wp $ and $ \wp^\prime $
		in equation \eqref{eq:wp-laurent}, and then determine the logarithm of a point
		belonging to $ \Imgof{\Exp} $ modulo $ p^k $.

		Specifically, by taking the series expansions in equation \eqref{eq:wp-laurent}
		modulo $ p^5 $, since $ p \mid z $, we have that
		\begin{align*}
			\wp(z) &\equiv \frac{1}{z^2} - \frac{1}{5} a z^2 - \frac{1}{7} b z^4 \modp{p^5},\\
			\frac{1}{2} \wp^\prime(z) &\equiv \frac{-1}{z^3} - \frac{1}{5} a z - \frac{2}{7} b z^3 \modp{p^5}.
		\end{align*}
		As $ \homcoords{Z}{X}{Y} = P \in \Imgof{\Exp} $ is such that
		$ \wp(z) \equiv \frac{X}{Z} \modp{p^k} $, and $ \frac{1}{2} \wp^\prime(z) \equiv \frac{Y}{Z} \modp{p^k} $,
		by multiplying both sides of the above equations by $ z^3 $, we obtain
		\begin{align*}
			z^3 \wp(z) &\equiv z \modp{p^5},\\
			\frac{z^3}{2} \wp^\prime(z) &\equiv -1 - \frac{1}{5} a z^4 \modp{p^5}.
		\end{align*}
		Hence, the first equivalence implies that $ z^2 \equiv \frac{Z}{X} \modp{p^{5}} $ since
		$ \wp(z) \equiv \frac{X}{Z} $, whereas the second equivalence give us
		\begin{align*}
			\frac{z^3}{2} \wp^\prime(z) &\equiv z^3 \frac{Y}{Z} \equiv z \frac{Z}{X} \frac{Y}{Z} = z \frac{Y}{X} \modp{p^{5}}, \\\nonumber
			\frac{z^3}{2} \wp^\prime(z) &\equiv - 1 - \frac{1}{5} a z^4 \equiv - 1 - \frac{1}{5} a {\roundbrack{\frac{Z}{X}}}^2 \modp{p^{5}},
		\end{align*}
		which implies
		\begin{equation*}
			z \equiv \frac{X}{Y} \roundbrack{-1 - \frac{1}{5} a {\roundbrack{\frac{Z}{X}}}^2} \modp{p^{5}}.
		\end{equation*}
		Note that if $ Z = 0 $ or $ k \le 4 $, then the logarithm of $ P \in \Imgof{\Exp} $ is
		$ - \frac{X}{Y} \modp{p^k} $.
		Since we determined the logarithm modulo $ p^5 $, we can compute the value of $ z $
		modulo $ p^k $ as $ z + h \cdot p^5 $, for $ h \in \N $, and $ k > 5 $.
		In order to determine the value of $ z $ modulo $ p^k $,
		since $ \wp(z) \equiv \frac{X}{Z} \modp{p^i} $,
		we compare the series expansion of $ \wp(z) $ modulo $ p^i $,
		for $ i = 6, \ldots, k $, with $ \frac{X}{Z} $. If this latter equivalence holds
		and $ \Expof{z} \ne P $, then $ z \modp{p^i} $ is \enquote{candidate} as a solution,
		and we can use it to compute the next candidates $ z^\prime = z + h^\prime \cdot p^{i + 1} $,
		for $ h^\prime = 1, \ldots, p - 1 $.
	\end{proof}

	\section{An attack to the ECDLP}\label{sec:ell-ecdlp}

	In this section, we show that the map $ \Log $ in \thmref{log-map} is closely related to the discrete logarithm
	problem for elliptic curves.

	\begin{theorem}[ECDLP over $ \Q $]\label{thm:ecdlp-q}
		Let $ \weier $ be an elliptic curve in Weierstrass form defined over field $ \Q $ of rational numbers,
		and let $ P $ and $ Q $ be two points belonging to $ \weierof{\Q} $ such that
		$ Q = h \cdot P $, for some $ h \in \N $. Lastly, let $ p $ be an odd prime number
		such that $ \weier $ has a good reduction modulo $ p $, and
		let $ t $ be equal to $ \ordp{\QModxof{}{P} - \Omega} $. If $ t > 2 $, then
		it is possible to determine the discrete logarithm of the point $ Q $ in at most
		$ \ceil{\log_p(h)} $ steps.
	\end{theorem}
	\begin{proof}
		Since $ Q - h \cdot P = \Omega $, then
		\begin{equation*}
			\QModxof{k}{Q - h \cdot P} \in \Imgof{\Exp}.
		\end{equation*}
		Since, over $ \Zof{p} $, the value of $ h $ is reduced modulo $ t $,
		we have that $ h = \overline{h} + n t $, for some $ n \in \N $, i.e. $ h \equiv \overline{h} \modp{t} $.
		Therefore, we have that $ Q = h \cdot P = (\overline{h} + nt) \cdot P $, and thus
		$ Q - \overline{h} \cdot P = nt \cdot P $.
		Since $ t \cdot \QModxof{k}{P} \in \Imgof{\Exp} $,
		we have that
		\begin{align*}
				\QModxof{k}{Q - \overline{h} \cdot P} &= \Expof{p^d l_2},\\
				n t \cdot \QModxof{k}{P} &= n \cdot \Expof{p^c l_1} = \Expof{p^c n l_1},\\
			\end{align*}
		for some $ l_1, l_2 \in \Z $, and $ 1 \le c, d \in \N $. Hence, by the taking the logarithm defined in
		\thmref{log-map} of both expressions, we have that
		\begin{equation*}
				\Logof{nt \cdot \QModxof{k}{P}} \equiv \Logof{\QModxof{k}{Q - \overline{h} \cdot P}} \modp{p^k},
			\end{equation*}
		which implies
		\begin{equation*}
				p^c n l_1 \equiv p^d l_2 \modp{p^k},
			\end{equation*}
		and, if $ l_1 \not\equiv 0 \modp{p^k} $, then
		\begin{align*}
				n \equiv p^{d - c} \frac{l_2}{l_1} \modp{p^{k - c - \nu_p(l_1)}} \quad \text{if } k - c - \nu_p(l_1) > 0,
			\end{align*}
		where $ \nu_p(\cdot) $ is the $ p $-adic evaluation function.

		We are left to compute the number of steps required to determine
		the correct value of $ n $ such that $ h = \overline{h} + n t $.
		Specifically, let $ n_2 \le n_3 \le \ldots \le n_\infty $ be the non decreasing sequence
		of the values of $ n $ computed modulo $ p^{k - c_k - \nu_p(l_{(1, k)})} $,
		for increasing values of $ 2 \le k $,
		where $ p^{c_k} l_{(1, k)} = \Logof{t \cdot \QModxof{k}{P}} $.
		Let $ i $ be the minimum index such that
		$ n_i = n_{i + 1} = n_{i + 2} = \ldots = n_\infty $. Note that this index
		exists and is such that $ h = \overline{h} + n_i t $, i.e. $ n = n_i $.
		Hence, the index $ i $ is lesser or equal to $ \ceil{\log_p(h)} $.
	\end{proof}

	\begin{remark}
		We remark that, in \thmref{ecdlp-q}, if $ t \mid h $,
		then $ \QModxof{}{Q} = \Omega $, which implies that $ \overline{h} = 0 $. Additionally,
		if $ \cardof{\qpjacobianof{\weier}{\Zof{p}}} $ is \enquote{small},
		then the value of $ \overline{h} $ can be easily determined.
	\end{remark}

	Hence, as a consequence of \thmref{ecdlp-q},
	in order to compute the value of $ n $ (and thus the value of $ h $),
	we could take simply the logarithms of $ t \cdot \QModxof{k}{P} $,
	and $ \QModxof{k}{Q - \overline{h} \cdot P} $, for increasing values of $ k $.

	Below, we give an example by using the above method.
	\begin{example}
		Let $ \weier $ be the elliptic curve in Weierstrass form
		over $ \Q $ given by the equation $ y^2 = x^3 - x + \frac{1}{4} $,
		let $ P = \roundbrack{2, \frac{5}{2}} $ be a point belonging to
		$ \weierof{\Q} $, and let $ Q = h \cdot P $. Suppose, for the purposes of this demonstration,
		we know that $ h = 31 $ and we want to verify the efficacy of the algorithm discussed above.

		Since the above curve has a good reduction modulo $ p = 3 $, we
		choose this small prime number. The reduction modulo $ 3 $ of the above
		curve is $ y^2 = x^3 - x + 1 $, whose Jacobian has cardinality equal to $ 7 $.
		Therefore, this group in cyclic and every point is a generator (in particular,
		$ \Modof{P} $).
		Since this group is small, it is easy to compute the reduction of
		$ h $ modulo $ t = \ordp{\Modof{P} - \Omega} = 7 $, and in particular we get $ \overline{h} = 3 $.
		Hence, by computing the logarithms of $ 7 \cdot \QModxof{k}{P} $,
		and $ \QModxof{k}{Q - 3 \cdot P} $, for increasing value
		of $ k $, we are able to compute the value of $ h $.
		In \tabref{ecdlp-ref}, we present the results of the aforementioned algorithm.
	\end{example}

	\begin{table}[tbh!]
		\centering
		\begin{tabular}{clZlZll}
			\hline\\
			$ k $ & $ \Logof{7 \cdot \QModof{3^k}{P}} $ & $ l_1 $ & $ \Logof{\QModof{3^k}{Q - 3 \cdot P} } $ & $ l_2 $ & $ n $ & $ h $ \\
			\hline\\
			1 & 3 & 1 & 0 & 0 & 0 & 3\\
			2 & 6 & 2 & 6 & 2 & 1 & 10\\
			3 & 6 & 2 & 24 & 8 & 4 & 31\\
		\end{tabular}
		\caption{ECDLP as successive approximations modulo $ 3^k $}
		\label{tab:ecdlp-ref}
	\end{table}

	\section{A generalized attack to the ECDLP}\label{sec:ell-dlp-Q-Qp}

	In this section, we provide a different way to use the results in \secref{ell-ecdlp}
	in order to compute the discrete logarithm for elliptic curves over $ \F_q $,
	where $ q $ is a big prime number.

	Let $ \weier $ be an elliptic curve in Weierstrass form over $ \Q $,
	and let $ P, Q $ be two points belonging to $ \weierof{\Q} $ such that
	$ Q = a \cdot P $, for some $ a \in \Z $.
	If $ \weier $ has a good reduction modulo $ q $, then the relationship
	between $ P $ and $ Q $ still holds over $ \F_q $, i.e.
	$ \QModof{q}{Q} = \overline{a} \cdot \QModof{q}{P} $,
	where $ a \equiv \overline{a} \modp{n} $, with $ n = \ordp{\QModof{q}{P} - \Omega} $.
	The same results hold if we consider a small prime number $ p \ll q $
	such that $ \weier $ has a good reduction modulo $ p $.
	The following diagram illustrates the concepts discussed previously:
	\begin{equation*}
		\begin{CD}
			\weierof{\Q} @>{\QMod_q}>> \weierof{\F_q} \\
			@VV{\QMod_p}V @. \\
			\weierof{\F_p} @. \\
		\end{CD}
	\end{equation*}
	where $ p $ and $ q $ are prime numbers such that $ q \gg p $.

	For instance, suppose one wants to compute the discrete logarithm
	of a point $ Q^\prime = a \cdot P^\prime $ belonging to the well-known
	Curve25519, where $ q = 2^{255} - 19 $.
	If we know an elliptic curve $ \weier $ over $ \Q $ such that
	the reduction of $ \weier $ modulo $ q $ is equal to Curve25519, there
	are two $ P, Q \in \weierof{\Q} $ such that $ \QModof{q}{P} = P^\prime $
	and $ \QModof{q}{Q} = Q^\prime $, and $ Q = b \cdot P $,
	where $ b \equiv a \modp{\ordp{P^\prime - \Omega}} $.
	Finding the discrete logarithm for Curve25519 is manifestly unpractical because
	of the cardinality of its Jacobian. However, since we supposed to have
	the curve $ \weier $, we may reduce it modulo a small prime $ p $
	(for instance, $ p = 3 $). Hence, we could compute the discrete logarithm
	of the point $ Q^\prime $ by using the method in \secref{ell-ecdlp}
	with the two points $ P $ and $ Q $. Specifically, we would find
	the integer $ b $ satisfying the equivalence $ Q = b \cdot P $,
	and we can compute the value $ a $ by reducing $ b $ modulo $ \ordp{P^\prime - \Omega} $.

	In practice, the elliptic curve over $ \Q $ and the
	relationship $ Q = b \cdot P $ are not typically available.
	One only has the elliptic curve over some finite field (such as Curve25519),
	a relationship between two points belonging to that curve, and wants to compute
	the discrete logarithm for that relationship.

	In principle, it would also be possible to derive the common ancestor of
	elliptic curve over $ \Q $, starting from the elliptic curve over $ \F_q $,
	but it is more complex.
	It turns out that if we have the relationship $ Q^\prime = a \cdot P^\prime $ over $ \F_q $,
	then this relationship over $ \Q $ is similar to
	$ Q = b \cdot P + \omega $, where $ b \equiv a \modp{\ordp{P^\prime - \Omega}} $, and
	$ \omega $ is a point belonging to $ \weierof{\Q} $
	which is equal to $ \Omega $ when reduced modulo $ q $.

	\section{DLP for hyperelliptic curves over \texorpdfstring{$\Qp{p}$}{Qp}}\label{sec:hec}

	In this section, we extend the method presented in \secref{ell-ecdlp} to
	hyperelliptic curves over $ \Qp{p} $, and we prove that these curves are resistant
	to this attack.

	\begin{definition}
		Let $ \K $ be a field. A hyperelliptic curve in Weierstrass form over $ \K $
		of genus $ g \ge 1 $ is an algebraic curve given by the equation
		\begin{equation}\label{eq:hec}
			y^2 + h(x) y = f(x),
		\end{equation}
		where $ \degof{h(x)} \le \floor{\frac{\degof{f(x)}}{2}} $,
		and $ f(x) $ is a monic polynomial such that
		$ \degof{f(x)} $ is either $ 2g + 1 $ or $ 2g + 2 $.
	\end{definition}

	We confine ourselves to the case in which $ \degof{f(x)} = 2 g + 1 $, i.e.
	the case of imaginary hyperelliptic curves.

	As for the elliptic curves, if $ \charof{\K} \ne 2, 3 $, then the map
	in equation \eqref{eq:tranform-map} transforms the equation \eqref{eq:hec} to the form
	\begin{equation}\label{eq:hec-reduced}
		y^2 = g(x),
	\end{equation}
	where $ \degof{g(x)} = 2 g + 1 $.

	In the following, we recall the Mumford representation for the elements of the Jacobian
	of a hyperelliptic curve.

	\begin{definition}[Mumford representation]
		Let $ \hyperec $ be a hyperelliptic curve, given by the equation
		$ \eqref{eq:hec} $, over an algebraically closed
		field $ \overline{\K} $. For any element $ D \in \jacobianof{\hyperec} $ there
		is a linearly equivalent reduced divisor $ D^\prime = \sum_{i = 1}^{n} P_i - n \cdot \Omega $,
		with $ 0 \le n \le g $, such that $ (x_i, y_i) = P_i \in \hyperecof{\overline{\K}} $,
		$ P_i \ne \Omega $, and $ P_i + P_j - 2 \cdot \Omega \not\equiv \Omega - \Omega $, for any $ i \ne j $.
		The reduced divisor $ D^\prime $ has a Mumford representation given by a pair of polynomials
		$ (u(x), v(x)) $ such that
		\begin{itemize}
			\itemsep0em
			\item $ \commandof{\gcd}{u(x), u^\prime(x), v(x)} = 1 $,
			\item $ \degof{v(x)} < \degof{u(x)} \le g $,
			\item $ u(x) \mid \roundbrack{v^2(x) + v(x) h(x) - f(x)} $,
		\end{itemize}
		where $ u(x) $ is a monic polynomial equal to $ \prod_{i = 1}^{n} (x - x_i) $,
		and $ v(x) $ is such that $ v(x_i) = y_i $ for all $ i $. Moreover,
		if $ (x_i, y_i) $ has multiplicity $ n_i > 1 $ in $ \suppof{D^\prime} $, then the following
		conditions must to be met
		\begin{equation*}
			\left.\roundbrack{{\roundbrack{\frac{\delta}{\delta x}}}^j \roundbrack{v^2(x) + v(x)h(x) - f(x)}}\right\rvert_{x = x_i} = 0
			\qquad \mbox{with } j = 0, \ldots, n_i - 1.
		\end{equation*}
	\end{definition}

	In the \exampleref{hec-fail}, we show that the following diagram is not commutative:
	\begin{equation}\label{eq:hec-diagram}
		\begin{CD}
			\qpjacobianof{\Q}{\hyperec} @. \times @. \qpjacobianof{\Q}{\hyperec} @>(+)>> \qpjacobianof{\Q}{\hyperec} \\
			@VV{\QModx{k}}V @. @VV{\QModx{k}}V @VV{\QModx{k}}V \\
			\qpjacobianof{\Zof{p^k}}{\hyperec} @. \times @. \qpjacobianof{\Zof{p^k}}{\hyperec} @>(+)>> \qpjacobianof{\Zof{p^k}}{\hyperec} \\
		\end{CD}
	\end{equation}
	where $ \qpjacobianof{\Q}{\hyperec} $ is the Jacobian of a hyperelliptic curve
	over $ \Q $. Hence, it is not possible to attack a hyperelliptic curve by using the
	method in \secref{ell-ecdlp}.

	\begin{example}\label{example:hec-fail}
		Let $ \hyperec $ be a hyperelliptic curve over $ \Q $ given by the equation
		\begin{equation*}
			y^2 = x(x - 1)(x + 2)(x - 3)(x + 4).
		\end{equation*}
		Additionally, let $ D $ be the reduced divisor in Mumford form equal to $ (x - 4, 24) $
		belonging to $ \qpjacobianof{\Q}{\hyperec} $. Lastly, let $ p $ be equal to $ 11 $.
		Note that, the curve $ \hyperec $ has a good reduction modulo $ p $ since its
		discriminant \cite{Lockhart} is different from zero modulo $ p $.

		Since the cardinality of $ \qpjacobianof{\Zof{11}}{\hyperec} $ is small,
		it is easy to compute the order of $ \QModxof{}{D} $, which is $ 16 $.

		Although on elliptic curves over $ \Q $ one has that
		$ h \cdot \QModxof{}{P} = \QModxof{}{h \cdot P} $, in the case of
		hyperelliptic curves the equivalent operations have two different results.
		On one hand, since $ \ordp{\QModxof{}{D}} = 16 $, then
		$ \QModxof{}{16 \cdot D} $ is equal to the neutral element
		$ \Omega - \Omega $, which is equal to $ (1, 0) $ in Mumford form.
		On the other hand, $ 16 \cdot \QModxof{}{D} $ is different from $ \QModxof{}{16 \cdot D} $.
		In particular, we have have that, from a certain $ h $ onwards, the
		divisors $ h \cdot \QModxof{}{D} $ and $ \QModxof{}{h \cdot D} $
		will be different. In this specific case, $ h $ is equal to $ 8 $
		and we have that
		\begin{align*}
			\QModxof{}{8 \cdot D} &= (x + 10, 4x + 7),\\
			8 \cdot \QModxof{}{D} &= (x + 10, 0).
		\end{align*}
		Note that the second result is the correct one since $ \QModxof{}{D} $ has order $ 16 $,
		and thus $ 8 \cdot \QModxof{}{D} $ has order $ 2 $, since $ v(x) = 0 $.
		Additionally, $ \QModxof{}{8 \cdot D} $ is not a valid divisor in Mumford
		form as $ \degof{4 x + 7} $ is not lesser than $ \degof{x + 10} $.
		Lastly, the divisor $ \QModxof{}{16 \cdot D} = (x^2, 8x) $ has not a
		valid Mumford form as well since $ \commandof{\gcd}{u(x), u^\prime(x), v(x)} \ne 1 $.
	\end{example}

	The source of the problem discussed in \exampleref{hec-fail} is the
	non-commutativity of $ \gcd $ and the reduction modulo $ p $, which
	are used in the Cantor-Koblitz algorithm \cite{cantor, Koblitz1989} for the addition
	of two divisors in Mumford form (see \algref{cantor-koblitz}) belonging to the Jacobian
	of a hyperelliptic curve.
	In particular, when $ p $ divides the coefficients of two polynomials
	$ f_1(x) $ or $ f_2(x) $, one has that
	\begin{equation*}
		\gcd(f_1(x), f_2(x)) \modp{p^k} \not\equiv \gcd(f_1(x) \modp{p^k}, f_2(x) \modp{p^k}).
	\end{equation*}

	Specifically, the Cantor-Koblitz algorithm computes the $ \gcd $ between the
	polynomials $ u_1(x) $ and $ u_2(x) $ of the Mumford form of the divisors to be summed, respectively.
	However, this computation may differ if one considers the two divisors over $ \Q $
	and their respective reductions modulo $ p^k $.

	In conclusion, the non-commutativity of the diagram in equation \eqref{eq:hec-diagram}
	implies that it is not possible to extend the results in \secref{ell-qp} to the case of
	hyperelliptic curves of genus $ g > 1 $, and as a consequence it is not possible
	to compute the DLP for hyperelliptic curves by using the method in \secref{ell-ecdlp}.

	\begin{algorithm}[tbh!]
		\label{alg:cantor-koblitz}
		\scriptsize
		\SetAlgoLined
		\SetNoFillComment
		\DontPrintSemicolon

		\SetKwData{ux}{$ u_1(x) $}
		\SetKwData{uux}{$ u_2(x) $}
		\SetKwData{dx}{$ d_1 $}
		\SetKwData{vvh}{$ v_1(x) + v_2(x) + h(x) $}
		\SetKwFunction{ExtGCD}{Extended-GCD}

		\KwIn{$ D_1 = (u_1(x), v_1(x)) $, $ D_2 = (u_2(x), v_2(x)) $, $ \hyperec: y^2 + h(x) y = f(x) $, where $ \deg(f(x)) = 2g + 1 $}
		\KwOut{$ D = D_1 + D_2 = (u(x), v(x)) $}
		\BlankLine
		\tcp{compute the extended GCD between $ u_1 (x) $ and $ u_2(x) $}
		$ (d_1, e_1, e_2) \longleftarrow \ExtGCD(\ux, \uux) $\tcp*{$ d_1 = e_1 u_1(x) + e_2 u_2(x) $}
		\tcp{compute the extended GCD between $ d_1 $ and $ v_1(x) + v_2(x) + h(x) $}
		$ (d, c_1, c_2) \longleftarrow \ExtGCD{\dx, \vvh} $\tcp*{$ d = c_1 d_1 + c_2 (v_1(x) + v_2(x) + h(x)) $}
		$ s_1 \longleftarrow c_1 e_1 $\;
		$ s_2 \longleftarrow c_1 e_2 $\;
		$ s_3 \longleftarrow c_2 $\;
		\tcp{composition step}
		$ u(x) \longleftarrow \frac{u_1(x) u_2(x)}{d^2} $\;
		$ v(x) \longleftarrow \frac{s_1 u_1(x) v_2(x) + s_2 u_2(x) v_1(x) + s_3 (v_1(x) v_2(x) + f(x))}{d} \modp{u(x)} $\;
		\Repeat(\tcp*[h]{reduction step}){$ \deg(u(x)) \le g $}{
			$ u(x) \longleftarrow \frac{f(x) - v(x) h(x) - v^2(x)}{u(x)} $\;
			$ v(x) \longleftarrow - h(x) - v(x) \modp{u(x)} $
		}
		$ u(x) \longleftarrow \frac{u(x)}{lc(u(x))} $\tcp*{divide $ u(x) $ by its leading coefficient}
		\KwRet{$ (u(x), v(x)) $}\;
		\caption{Cantor-Koblitz algorithm}
	\end{algorithm}

	\section{Conclusions}\label{sec:conclusion}

	In this paper, we extended the exponential map $ \Exp $ for elliptic curves
	over $ \Zof{p^k} $. Additionally, we gave an explicit method in \thmref{log-map}
	to compute the inverse map
	of $ \Exp $, denoted as $ \Log $, and an algorithm in \thmref{ecdlp-q}
	to attack the elliptic curve discrete logarithm by using the maps
	$ \Exp $ and $ \Log $ over $ \Zof{p^k} $.
	Lastly, we proved that hyperelliptic curves are resistant to our method, meaning that
	these curves provide a higher security compared to that of the elliptic curves.

	\printbibliography[title={References}]

\end{document}